\documentclass{amsart}

\newtheorem{theorem}{Theorem}[section]
\newtheorem{lemma}[theorem]{Lemma}

\theoremstyle{definition}

\theoremstyle{remark}
\newtheorem{remark}[theorem]{Remark}

\numberwithin{equation}{section}

\newcommand{\tmmathbf}[1]{\ensuremath{\boldsymbol{#1}}}

\begin{document}

\title{On the solvability of systems of bilinear equations in finite
fields}

\author{Le Anh Vinh}
\address{Department of Mathematics, Harvard University, Cambridge, MA 02138}
\email{vinh@math.harvard.edu}



\date{Dec 1, 2008}



\begin{abstract}
Given $k$ sets $\mathcal{A}_i \subseteq \mathbb{F}_q^d$ and a
  non-degenerate bilinear form $B$ in $\mathbb{F}_q^d$. We consider the system of
  $l \leq \binom{k}{2}$ bilinear equations
  \[ B (\tmmathbf{a}_i, \tmmathbf{a}_j) = \lambda_{i j},\, \, \, \tmmathbf{a}_i \in
     \mathcal{A}_i, i = 1, \ldots, k. \]
  We show that the system is solvable for any $\lambda_{i j} \in \mathbb{F}_q^{*}$, $1 \leq i,j \leq k$, given that the restricted sets
  $\mathcal{A}_i$'s are sufficiently large.
\end{abstract}

\maketitle

\section{Introduction}

In \cite{s1}, S\'ark\"ozy proved that if $\mathcal{A}$, $\mathcal{B}$, $\mathcal{C}$, $\mathcal{D}$ are ``large'' subsets of $\mathbb{Z}_p$, more precisely, $|\mathcal{A}| |\mathcal{B}|  |\mathcal{C}| |\mathcal{D}| \gg p^3$, then the equation 
\begin{equation}\label{e1}
ab + 1 = cd
\end{equation}
can be solved with $a \in \mathcal{A}$, $b\in \mathcal{B}$, $c\in \mathcal{C}$ and $d \in \mathcal{D}$. Gyarmati and S\'ark\"ozy \cite{s2} generalized the results on the solvability of equation (\ref{e1}) to finite fields. They also study the solvability of other (higher degree) algebraic equations with solutions restricted to ``large'' subsets of $\mathbb{F}_q$, where $\mathbb{F}_q$ denote the finite field of $q$ elements. Using exponential
sums, Hart and Iosevich \cite{hart-iosevich} studied similar problem for any bilinear equation over $\mathbb{F}_q^d$. They showed that for any two sets $\mathcal{A},
\mathcal{B} \subseteq \mathbb{F}_q^d$, with
\[ |\mathcal{A}| |\mathcal{B}| > C q^{d + 1} \]
for some absolute constant $C > 0$, the equation
\[ \tmmathbf{a} \cdot \tmmathbf{b}= \lambda,\,\,\, \tmmathbf{a} \in \mathcal{A},
   \tmmathbf{b} \in \mathcal{B}, \]
is solvable for any $\lambda \in \mathbb{F}_q^{\ast}$. Although the proof is
given only in the case of the dot product, it goes through without any
essential changes if the dot product $\tmmathbf{a} \cdot \tmmathbf{b}$ is
replaced by any non-degenerate bilinear form $B (\tmmathbf{a}, \tmmathbf{b})$. Using bounds of multiplicative character sums, Shparlinski \cite{igor} extended the class of sets which satisfy this property. 

In this paper, we will use methods from graph theory to study the system of bilinear equations in finite
fields. More precisely, we consider the following system
\[ B (\tmmathbf{a}_i, \tmmathbf{a}_j) = \lambda_{i j}, \,\,\, \tmmathbf{a}_i \in
   \mathcal{A}_i, i = 1, \ldots, k \]
over $\mathbb{F}_q^d$, with variables from arbitrary sets $\mathcal{A}_i
\subseteq \mathbb{F}_q^d$, $i = 1, \ldots, k$. Our first result is the
following.

\begin{theorem}\label{general}
  Given $k$ sets $\mathcal{A}_i \subseteq \mathbb{F}_q^d$ and let $B (\cdot,
  \cdot)$ be a non-degenerate bilinear form in $\mathbb{F}_q^d$. Consider the
  system $\mathcal{L}$ of $l \leq \binom{k}{2}$ bilinear equations
  \begin{equation}\label{gs} B (\tmmathbf{a}_i, \tmmathbf{a}_j) = \lambda_{i j}, \,\,\, \tmmathbf{a}_i \in
     \mathcal{A}_i, i = 1, \ldots, k. \end{equation}
  Suppose that each variable appears in at most $t \leq k - 1$ equations
  and
  \[ |\mathcal{A}_i | \gg q^{\frac{d - 1}{2} + t} . \]
  Then for any $\lambda_{i j} \in \mathbb{F}_q^{\ast}$, the above system has
  \[ (1 + o (1)) q^{- l} \prod_{i = 1}^k |\mathcal{A}_i | \]
  solutions. 
\end{theorem}

The study of systems of bilinear equations in vector spaces over finite fields in the context of Theorem \ref{general} in the case $t=k-1$ can be found in Chapter 2 of D. Hart's dissertation \cite{hart-thesis}. Theorem \ref{general} is a quantitative improvement over the result given there.

The serious difficulties arise when the number of equations that each
variable involves is sufficiently large with respect to the ambient dimension.
More precisely, Theorem \ref{general} is only non-trivial in the range $d \geq 2
t$. Even in the case of each variable appears in at most two equations,
Theorem \ref{general} is non non-trivial for $d \geq 4$. The purpose of the
following theorem is to fill in this gap. We show that the
system of two bilinear equations and three variables in large restricted
subsets of $\mathbb{F}_q^d$ is always solvable for $d \geq 2$.

\begin{theorem} \label{23}
  Given three sets $\mathcal{A}, \mathcal{B}, \mathcal{C} \subseteq
  \mathbb{F}_q^d$ and let $B (\cdot, \cdot)$ be a non-degenerate bilinear
  form in $\mathbb{F}_q^d$. Consider the system of two equations
  \begin{equation}\label{23s} B (\tmmathbf{a}, \tmmathbf{b}) = \lambda_1, B (\tmmathbf{a}, \tmmathbf{c}) =
     \lambda_2, \,\,\, \tmmathbf{a} \in \mathcal{A}, \tmmathbf{b} \in \mathcal{B},
     \tmmathbf{c} \in \mathcal{C}. \end{equation}
  Suppose that
  \[ |\mathcal{A}| |\mathcal{B}|, |\mathcal{A}| |\mathcal{C}| \gg q^{d + 1},
  \]
  then for any $\lambda_1, \lambda_2 \in \mathbb{F}_q^{\ast}$, the above
  system has
  \[ (1 + o (1)) \frac{|\mathcal{A}| |\mathcal{B}| |\mathcal{C}|}{q^2} \]
  solutions.
\end{theorem}

The solvability of the system of three bilinear equations and three variables in large restricted
subsets of $\mathbb{F}_q^2$, however, is harder to determine. We will instead
show that the system is solvable for a positive proportion of all triples
$(\lambda_1, \lambda_2, \lambda_3) \in \mathbb{F}_q^{\ast}$. More precisely,
we have the following theorem.

\begin{theorem} \label{332}
  Given three sets $\mathcal{A}, \mathcal{B}, \mathcal{C} \subseteq
  \mathbb{F}_q^2$ and let $B (\cdot, \cdot)$ be a non-degenerate bilinear
  form in $\mathbb{F}_q^2$. Consider the system of equations
  \[ B (\tmmathbf{a}, \tmmathbf{b}) = \lambda_1, B (\tmmathbf{a}, \tmmathbf{c}) =
     \lambda_2, B (\tmmathbf{b}, \tmmathbf{c}) = \lambda_3,\,\,\, \tmmathbf{a} \in
     \mathcal{A}, \tmmathbf{b} \in \mathcal{B}, \tmmathbf{c} \in \mathcal{C}. \]
  Suppose that
  \[ |\mathcal{A}|, |\mathcal{B}|, |\mathcal{C}| \gg q^{3 / 2}, \]
  then the above system is solvable for $\Omega ( \frac{\sqrt{|\mathcal{B}|
  |\mathcal{C}|}}{q^2}) q^3$ \ triples $(\lambda_1, \lambda_2, \lambda_3) \in
  (\mathbb{F}_q^{\ast})^3$. \ 
\end{theorem}

It is conceivable that we can chop off the term $\Omega (
\frac{\sqrt{|\mathcal{B}| |\mathcal{C}|}}{q^2})$ in the above theorem, or even
better, the system is solvable for $(1 - o (1)) q^3$ triples
$(\lambda_1, \lambda_2, \lambda_3) \in (\mathbb{F}_q^{\ast})^3$. We show that
it is indeed the case when the ambient space has the dimension $d \geq
3$.

\begin{theorem} \label{333}
  Given three sets $\mathcal{A}, \mathcal{B}, \mathcal{C} \subseteq
  \mathbb{F}_q^d$ and let $B (\cdot, \cdot)$ be a non-degenerate bilinear
  form in $\mathbb{F}_q^d$. Consider the system of equations
  \[ B(\tmmathbf{a},\tmmathbf{b})= \lambda_1, B(\tmmathbf{a},\tmmathbf{c})=
     \lambda_2, B(\tmmathbf{b},\tmmathbf{c})= \lambda_3, \,\, \tmmathbf{a} \in
     \mathcal{A}, \tmmathbf{b} \in \mathcal{B}, \tmmathbf{c} \in \mathcal{C}. \]
  Suppose that
  \[ |\mathcal{A}| |\mathcal{B}|, |\mathcal{A}| |\mathcal{C}|, |\mathcal{B}|
     |\mathcal{C}| \gg q^{d + 2}, \]
  then the above system is solvable for $(1 - o (1)) q^3$ triples $(\lambda_1,
  \lambda_2, \lambda_3) \in (\mathbb{F}_q^{\ast})^3$.
\end{theorem}

Note that there is a series of papers dealing with similar results in the sovability of systems of quadratic forms, for example, see \cite{covert,hart2,hart3,vinh1,vinh2,vinh3,vinh4}.
\section{Bilinear equations in finite fields}

Let $B (\cdot, \cdot)$ be a nondegenerate bilinear form in $\mathbb{F}_q^d$
and $\lambda \in \mathbb{F}_q^{\ast}$. For any $\tmmathbf{v} \in
\mathbb{F}_q^d$ and a subset $V \subseteq \mathbb{F}_q^d$, denote
$N^{\lambda} (\tmmathbf{v})$ be the set of all vectors $\tmmathbf{u} \in
\mathbb{F}_q^d$ such that $B (\tmmathbf{v}, \tmmathbf{u}) = \lambda$, and let
$N^{\lambda}_V (\tmmathbf{v}) = N^{\lambda} (\tmmathbf{v}) \cap V$. The
following key estimate says that the cardinalities of $N^{\lambda}_V
(\tmmathbf{v})$'s are close to $|V| / q$ when $|V|$ is large.

\begin{lemma} \label{key}
  For every subset $V$ of $\mathbb{F}_q^d$ then
  \[ \sum_{\tmmathbf{v} \in \mathbb{F}_q^d} \left( |N_V^{\lambda} (\tmmathbf{v}) | -
     \frac{|V|}{q} \right)^2 < q^{d - 1} |V|. \]
\end{lemma}

\begin{proof}
  For any set $X$, let $X (\cdot)$ denote the characteristic function of $X$.
  Let $\chi$ be any non-trivial additive character of $\mathbb{F}_q$. We have
  \begin{eqnarray*}
    |N_V^{\lambda} (\tmmathbf{v}) | & = & \sum_{\tmmathbf{u} \in
    \mathbb{F}_q^d, B (\tmmathbf{v}, \tmmathbf{u}) = \lambda} V (\tmmathbf{u})\\
    & = & \sum_{u \in \mathbb{F}_q^d, s \in \mathbb{F}_q} \frac{1}{q} \chi
    (s (B (\tmmathbf{v}, \tmmathbf{u}) - \lambda)) V (\tmmathbf{u})\\
    & = & \frac{|V|}{q} + \frac{1}{q} \sum_{u \in \mathbb{F}_q^d, s \in
    \mathbb{F}_q^{\ast}} \chi (s (B (\tmmathbf{v}, \tmmathbf{u}) - \lambda)) V
    (\tmmathbf{u}) .
  \end{eqnarray*}
  Therefore
  \begin{eqnarray}
    & & \sum_{\tmmathbf{v} \in \mathbb{F}_q^d} \left( |N_V^{\lambda} (\tmmathbf{v}) |  - 
    \frac{|V|}{q} \right)^2  =  \frac{1}{q^2} \sum_{\tmmathbf{v} \in
    \mathbb{F}_q^d} \left( \sum_{\tmmathbf{u} \in \mathbb{F}_q^d, s \in
    \mathbb{F}_q^{\ast}} \chi (s (B (\tmmathbf{v}, \tmmathbf{u}) - \lambda)) V
    (\tmmathbf{u}) \right)^2 \nonumber\\
    & = & \frac{1}{q^2} \sum_{\tmmathbf{v}, \tmmathbf{u}, \tmmathbf{u}' \in
    \mathbb{F}_q^d, s, s' \in \mathbb{F}_q^{\ast}} \chi ((s B (\tmmathbf{v},
    \tmmathbf{u}) - s' B (\tmmathbf{v}, \tmmathbf{u}')) \chi (\lambda (s' - s)) V
    (\tmmathbf{u}) V (\tmmathbf{u}') \nonumber\\
    & = & \frac{1}{q^2} (R_1 + R_2), \label{h1}
  \end{eqnarray}
  where $R_1$ is taken over $s = s'$ and $R_2$ is taken over $s \neq s'$. We
  compute each term.
  \begin{eqnarray}
    R_1 & = & \sum_{\tmmathbf{v}, \tmmathbf{u}, \tmmathbf{u}' \in
    \mathbb{F}_q^d, s = s' \in \mathbb{F}_q^{\ast}} \chi ((s B (\tmmathbf{v},
    \tmmathbf{u}) - s' B (\tmmathbf{v}, \tmmathbf{u}')) \chi (\lambda (s' - s)) V
    (\tmmathbf{u}) V (\tmmathbf{u}')\nonumber\\
    & = & \sum_{\tmmathbf{v}, \tmmathbf{u}, \tmmathbf{u}' \in \mathbb{F}_q^d, s
    = s' \in \mathbb{F}_q^{\ast}} \chi (s B (\tmmathbf{v},
    \tmmathbf{u}-\tmmathbf{u}')) V (\tmmathbf{u}) V (\tmmathbf{u}')\nonumber\\
    & = & (q - 1) q^d |V|,\label{h2}
  \end{eqnarray}
  where the last line follows from the orthogonality in $\tmmathbf{v}$. Now we
  compute $R_2$.
  \begin{eqnarray}
    R_2 & = & \sum_{\tmmathbf{v}, \tmmathbf{u}, \tmmathbf{u}' \in
    \mathbb{F}_q^d, s \in \mathbb{F}_q^{\ast}, a \neq 0, 1} \chi (s B
    (\tmmathbf{v}, \tmmathbf{u}- a\tmmathbf{u}')) \chi (\lambda (a s - s)) V
    (\tmmathbf{u}) V (\tmmathbf{u}')\nonumber\\
    & = & \sum_{\tmmathbf{v}, \tmmathbf{u}, \tmmathbf{u}' \in \mathbb{F}_q^d, s
    \in \mathbb{F}_q^{\ast}, a \neq 0, 1, \tmmathbf{u}= a\tmmathbf{u}'} \chi
    (\lambda (a s - s)) V (\tmmathbf{u}) V (a^{- 1} \tmmathbf{u})\nonumber\\
    & = & - \sum_{\tmmathbf{v}, \tmmathbf{u} \in \mathbb{F}_q^d, a \neq 0, 1}
    V (\tmmathbf{u}) V (a^{- 1} \tmmathbf{u})\nonumber\\
    & \geq & - (q - 2) q^d |V|,\label{h3}
  \end{eqnarray}
  where the second line follows from the orthogonality in $\tmmathbf{v}$. The
  lemma follows immediately from (\ref{h1}), (\ref{h2}) and (\ref{h3}). 
\end{proof}

The following result (\cite[Theorem 2.1]{hart-iosevich}) is an easy corollary of Lemma \ref{key}.

\begin{theorem} \label{12}
  (\cite[Theorem 2.1]{hart-iosevich}) Let $B (\cdot, \cdot)$ be a nondegenerate bilinear
  form in $\mathbb{F}_q^d$ and $\lambda \in \mathbb{F}_q^{\ast}$. For any
  two subsets $V, U \subseteq \mathbb{F}_q^d$, denote $N^{\lambda} (V, U)$ be
  the set of pairs $(\tmmathbf{v}, \tmmathbf{u}) \in V \times U$ such that $B
  (\tmmathbf{v}, \tmmathbf{u}) = \lambda$. Then we have
  \[ \left| N^{\lambda} (V, U) - \frac{|V| |U|}{q} \right| < \sqrt{q^{d - 1}
     |V| |U|} . \]
\end{theorem}

\begin{proof}
  By Lemma \ref{key}, we have
  \[ \sum_{\tmmathbf{u} \in U} \left( |N_V^{\lambda} (\tmmathbf{u}) | - \frac{|V|}{q}
     \right)^2 \leq \sum_{\tmmathbf{u} \in \mathbb{F}_q^d} \left( |N_V^{\lambda}
     (\tmmathbf{u}) | - \frac{|V|}{q} \right)^2 < q^{d - 1} |V|. \]
  By the Cauchy-Schwartz inequality,
  \begin{eqnarray*}
    \left| N^{\lambda} (V, U) - \frac{|V| |U|}{q} \right| & \leq &
    \sum_{\tmmathbf{u} \in U} \left| |N_V^{\lambda} (\tmmathbf{u}) | - \frac{|V|}{q}
    \right|\\
    & \leq & \sqrt{|U|} \sqrt{\sum_{\tmmathbf{u} \in U} \left( |N_V^{\lambda}
    (\tmmathbf{u}) - \frac{|V|}{q} \right)^2}\\
    & \leq & \sqrt{q^{d - 1} |V| |U|} .
  \end{eqnarray*}
\end{proof}

\begin{remark} Theorem \ref{12} has several applications in additive combinatorics (see \cite{hart}). We present here another application of this theorem to Waring's problem (mod $p$). Let $p$ be a prime and $k$ a positive integer. The smallest $s$ such that the congruence
\begin{equation} \label{w}
x_1^k + \ldots + x_s^k \equiv a \mod p
\end{equation} 
is solvable for all numbers $a$ is called Waring's number (mod $p$), denote $\gamma(k,p)$ (see \cite{cipra} for a historical background and recent results on this problem). Let $\gamma^{*}(k,p)$ denote the smallest $s$ such that the congruence (\ref{w}) is solvable for all $a \neq 0$. It is clear that $\gamma(k,p) \leq \gamma^{*}(k,p)+1$. Take $A$ be the set of all $k$-power in $\mathbb{F}_q$, and let $U \equiv V \equiv A^d$. From Theorem \ref{12}, the congruence (\ref{w}) with $s=d$ is solvable for all $a \neq 0$ if $q^{\frac{d-1}{2d}} \geq k$. It follows that
$\gamma^{*}(k,p) \leq d$ whenever $\gamma^{*}(k,p)$. This essentially matches with the classical bound of Weil in \cite{weil}.
  
\end{remark}

\section{General systems of bilinear equations}

  Now we give a proof of Theorem \ref{general}. The proof is very similar to that of \cite[Theorem~4.10]{pseudo}. Consider a random one-to-one mapping of the set of
  variables of $\mathcal{L}$ into the sets $\mathcal{A}_1,\ldots,\mathcal{A}_k$. Let $M (\mathcal{L})$ denote the event
  that all equations in $\mathcal{L}$ are satisfied under this mapping. We say that the mapping is an embedding of $\mathcal{L}$ in such a case. It
  suffices to prove that
  \begin{equation} \label{5}
    \Pr (M (\mathcal{L})) = (1 + o (1)) q^{-l}.
  \end{equation}
  We prove (\ref{5}) by induction on $l$, the number of equations in $\mathcal{L}$. The base case $l =
  0$ is trivial. Suppose that the theorem holds for all systems of less
  than $l$ equations. Let $B(\tmmathbf{a}_1,\tmmathbf{a}_2)=\lambda$, $\tmmathbf{a}_1 \in \mathcal{A}_1$, $\tmmathbf{a}_2 \in \mathcal{A}_2$, be an equation of $\mathcal{L}$. Let $\mathcal{L}_{\tmmathbf{a}_1}, \mathcal{L}_{\tmmathbf{a}_2}, \mathcal{L}_{ \{ \tmmathbf{a}_1, \tmmathbf{a}_2 \} }$ be the subsystems obtained from $\mathcal{L}$ by removing all equations that contain $\tmmathbf{a}_1$, $\tmmathbf{a}_2$, $\{\tmmathbf{a}_1, \tmmathbf{a}_2\}$, respectively. And let $\mathcal{L}_{\tmmathbf{a}_1,\tmmathbf{a}_2}$ be the subsystem obtained
  from $\mathcal{L}$ by removing the equation $B(\tmmathbf{a}_1,\tmmathbf{a}_2)=\lambda$. We have
  \begin{equation}
    \text{$\Pr (M (\mathcal{L}_{\tmmathbf{a}_1,\tmmathbf{a}_2})) = \Pr$} (M (\mathcal{L}_{\tmmathbf{a}_1,\tmmathbf{a}_2}) |M (\mathcal{L}_{ \{ \tmmathbf{a}_1, \tmmathbf{a}_2 \} })) \cdot \Pr (M
    (\mathcal{L}_{ \{ \tmmathbf{a}_1, \tmmathbf{a}_2 \} })) .
  \end{equation}
  Let $l_1$ be the number of equations in $\mathcal{L}_{ \{ \tmmathbf{a}_1, \tmmathbf{a}_2 \} }$. Since (\ref{5}) holds for
  $\mathcal{L}_{\tmmathbf{a}_1,\tmmathbf{a}_2}$ and $\mathcal{L}_{ \{ \tmmathbf{a}_1, \tmmathbf{a}_2 \} }$, we have \[ \Pr (M (\mathcal{L}_{\tmmathbf{a}_1,\tmmathbf{a}_2})) = (1 + o (1))
  q^{1-l},\] and \[ \Pr (M (\mathcal{L}_{ \{ \tmmathbf{a}_1, \tmmathbf{a}_2 \} })) = (1 + o (1))q^{-l_1}.\] Therefore, we have
  \begin{equation}
    \Pr (M (\mathcal{L}_{\tmmathbf{a}_1,\tmmathbf{a}_2}) |M (\mathcal{L}_{ \{ \tmmathbf{a}_1, \tmmathbf{a}_2 \} })) = (1 + o (1)) q^{l_1 + 1 - l} .
  \end{equation}
  For an embedding $f_1$ of $\mathcal{L}_{ \{ \tmmathbf{a}_1, \tmmathbf{a}_2 \} }$ into $\mathcal{A}_1,\ldots,\mathcal{A}_k$, let $\phi (\tmmathbf{a}_1, f_1)$,
  $\phi (\tmmathbf{a}_2, f_1)$ and $\phi (\tmmathbf{a}_1 \tmmathbf{a}_2, f_1)$ be the number of extensions of $f_1$ to an embedding of
  $\mathcal{L}_{\tmmathbf{a}_1}, \mathcal{L}_{\tmmathbf{a}_2}$ and $\mathcal{L}_{\tmmathbf{a}_1,\tmmathbf{a}_2}$ into $\mathcal{A}_1,\ldots,\mathcal{A}_k$, respectively. Note that an
  extension $f_{\tmmathbf{a}_1}$ of $f_1$ to an embedding of $\mathcal{L}_{\tmmathbf{a}_1}$ and an extension $f_{\tmmathbf{a}_2}$ of
  $f_1$ to an embedding of $\mathcal{L}_{\tmmathbf{a}_2}$ give us a unique extension of $f_1$ to an
  embedding of $\mathcal{L}_{\tmmathbf{a}_1,\tmmathbf{a}_2}$. Hence,
  \begin{equation}
    \phi (\tmmathbf{a}_1 \tmmathbf{a}_2, f_1) = \phi (\tmmathbf{a}_1, f_1)\phi (\tmmathbf{a}_2, f_1).
  \end{equation}
  Averaging over all possible extensions of $f_1$ to a mapping from $\mathcal{L}_{\tmmathbf{a}_1,\tmmathbf{a}_2}$ into $\mathcal{A}_1,\ldots,\mathcal{A}_k$, we have
  \begin{equation*}
     \Pr (M (\mathcal{L}_{\tmmathbf{a}_1,\tmmathbf{a}_2}) |f_1) =
    \frac{\phi (\tmmathbf{a}_1, f_1)\phi (\tmmathbf{a}_2, f_1)}{|\mathcal{A}_1| |\mathcal{A}_2|}.
  \end{equation*}
  Taking expectation over all embedding $f_1$, the LHS becomes
  \[
    \Pr (M (\mathcal{L}_{\tmmathbf{a}_1,\tmmathbf{a}_2}) |M (\mathcal{L}_{ \{ \tmmathbf{a}_1, \tmmathbf{a}_2 \} })) = (1 + o (1)) q^{l_1 + 1 - l} .
  \]
  So we get
  \begin{equation}\label{v4}
    E_{f_1} (\phi (\tmmathbf{a}_1, f_1) \phi (\tmmathbf{a}_2, f_1) |M (\mathcal{L}_{ \{ \tmmathbf{a}_1, \tmmathbf{a}_2 \} })) = (1 + o (1)) |\mathcal{A}_1| |\mathcal{A}_2|
    q^{l_1 + 1 - l} .
  \end{equation}
  Now, let $f$ be a random one-to-one mapping of the set of
  variables of $\mathcal{L}$ into the sets $\mathcal{A}_1,\ldots,\mathcal{A}_k$. Let $f_1$
  be a fixed embedding of $\mathcal{L}_{ \{ \tmmathbf{a}_1, \tmmathbf{a}_2 \} }$. Let $\mathcal{A}'_2$ and $\mathcal{A}'_1$ be the set of
  all possible images of $\tmmathbf{a}_2$ and $\tmmathbf{a}_1$ over all possible extensions of $f_1$ to
  embeddings of $\mathcal{L}_{\tmmathbf{a}_1}$ and $\mathcal{L}_{\tmmathbf{a}_2}$ into $\mathcal{A}_1,\ldots,\mathcal{A}_k$, respectively. From Theorem \ref{12}, the number of possible pairs $(\tmmathbf{a}_1, \tmmathbf{a}_2)$
  with $\tmmathbf{a}_1 \in \mathcal{A}'_1$ and $\tmmathbf{a}_2 \in \mathcal{A}'_2$ such that $B(\tmmathbf{a}_1,\tmmathbf{a}_2) = \lambda$ is bounded
  by
  \begin{equation} \label{11}
    \frac{\phi (\tmmathbf{a}_1, f_1)\phi (\tmmathbf{a}_2, f_1)}{q} \pm \sqrt{q^{d-1}\phi (\tmmathbf{a}_1, f_1)\phi (\tmmathbf{a}_2, f_1)}
  \end{equation}
  Thus, we have
  \[
    \Pr\ _f (M (\mathcal{L}) |f_{|\{\tmmathbf{a}_3,\ldots, \tmmathbf{a}_k\}} = f_1) = \frac{\phi (\tmmathbf{a}_1, f_1)\phi (\tmmathbf{a}_2, f_1)}{q|\mathcal{A}_1| |\mathcal{A}_2|} + \delta,
  \]
  where \[| \delta | \leq  \frac{\sqrt{q^{d-1}\phi (\tmmathbf{a}_1, f_1)\phi (\tmmathbf{a}_2, f_1)}}{|\mathcal{A}_1| |\mathcal{A}_2|}.\] Averaging over all possible embeddings
  $f_1$, we get
  \begin{eqnarray*}
    \Pr (M (\mathcal{L}) |M (\mathcal{L}_{ \{ \tmmathbf{a}_1, \tmmathbf{a}_2 \} }) & =& \frac{E_{f_1} (\phi (\tmmathbf{a}_1, f_1) \phi (\tmmathbf{a}_2, f_1) |M (\mathcal{L}_{ \{ \tmmathbf{a}_1, \tmmathbf{a}_2 \} }))}{q|\mathcal{A}_1| |\mathcal{A}_2|} + E_{f_1} (\delta)\\
    & = & (1 + o (1))q^{l_1 - l} + E_{f_1}
    (\delta),
  \end{eqnarray*}
  where the second lines follows from (\ref{v4}) and (\ref{11}). By Jensen's inequality, we have
  \begin{equation}
    |E_{f_1} (\delta) | \leq q^{(d-1)/2} \frac{\sqrt{E (\phi (\tmmathbf{a}_1, f_1)\phi (\tmmathbf{a}_2, f_1))}}{|\mathcal{A}_1| |\mathcal{A}_2|} = (1 + o (1)) \frac{q^{(d-1)/2}}{\sqrt{|\mathcal{A}_1| |\mathcal{A}_2|}}
    q^{(l_1 + 1 - l) / 2},
  \end{equation}
  which is negligible to the first term as \[ \sqrt{|\mathcal{A}_1| |\mathcal{A}_2|} \gg q^{\frac{d-1}{2}+t} \geq q^{(d-1)/2} q^{(l_1 + 1 - l)
  / 2}.\]
  Thus, we have
  \begin{equation}
    \Pr (M (\mathcal{L})) = \Pr (M (\mathcal{L}) |M (\mathcal{L}_{ \{ \tmmathbf{a}_1, \tmmathbf{a}_2 \} })  \Pr (M (\mathcal{L}_{ \{ \tmmathbf{a}_1, \tmmathbf{a}_2 \} }) = (1 + o
    (1)) q^{-l}.
  \end{equation}
  This completes the proof of the theorem.

\section{The system of two equations and three variables }

We will prove Theorem \ref{23} in this section. Our proof relies on Lemma \ref{key}
above. For any $\tmmathbf{v} \in \mathbb{F}_q^d$ and a subset $V \subseteq
\mathbb{F}_q^d$, denote $N^{\lambda} (\tmmathbf{v})$ be the set of all vectors
$\tmmathbf{u} \in \mathbb{F}_q^d$ such that $B (\tmmathbf{v}, \tmmathbf{u}) =
\lambda$, and let $N^{\lambda}_V (\tmmathbf{v}) = N^{\lambda} (\tmmathbf{v})
\cap V$. The number of solutions of the system (\ref{23s}) is
\[ \sum_{\tmmathbf{a} \in \mathcal{A}} |N_{\mathcal{B}}^{\lambda_1} (a) |
   |N_{\mathcal{C}}^{\lambda_2} (a) |. \]
From Lemma \ref{key}, we have
\begin{eqnarray*}
  \sum_{\tmmathbf{a} \in \mathcal{A}} \left( |N_{\mathcal{B}}^{\lambda_1}
  (\tmmathbf{a}) | - \frac{|\mathcal{B}|}{q} \right)^2 & \leq &
  \sum_{\tmmathbf{a} \in \mathbb{F}_q^d} \left( |N_{\mathcal{B}}^{\lambda_1}
  (\tmmathbf{a}) | - \frac{|\mathcal{B}|}{q} \right)^2 < q^{d - 1}
  |\mathcal{B}|\\
  \sum_{\tmmathbf{a} \in \mathcal{A}} \left( |N_{\mathcal{C}}^{\lambda_2}
  (\tmmathbf{a}) | - \frac{|\mathcal{C}|}{q} \right)^2 & \leq &
  \sum_{\tmmathbf{a} \in \mathbb{F}_q^d} \left( |N_{\mathcal{C}}^{\lambda_2}
  (\tmmathbf{a}) | - \frac{|\mathcal{C}|}{q} \right)^2 < q^{d - 1}
  |\mathcal{C}|.
\end{eqnarray*}
Thus, by the Cauchy-Schwartz inequality, we have
\begin{eqnarray*}
  & & \left[ \sum_{\tmmathbf{a} \in \mathcal{A}} \left(
  |N_{\mathcal{B}}^{\lambda_1} (\tmmathbf{a}) | - \frac{|\mathcal{B}|}{q}
  \right) \left( |N_{\mathcal{C}}^{\lambda_2} (\tmmathbf{a}) | -
  \frac{|\mathcal{C}|}{q} \right) \right]^2 \\
  &  & \leq  \sum_{\tmmathbf{a}
  \in \mathcal{A}} \left( |N_{\mathcal{B}}^{\lambda_1} (\tmmathbf{a}) | -
  \frac{|\mathcal{B}|}{q} \right)^2 \sum_{\tmmathbf{a} \in \mathcal{A}} \left(
  |N_{\mathcal{C}}^{\lambda_2} (\tmmathbf{a}) | - \frac{|\mathcal{C}|}{q}
  \right)^2 <  q^{2 d - 2} |\mathcal{B}| |\mathcal{C}|.
\end{eqnarray*}
This implies that
\begin{equation}\label{v1} \left| \sum_{\tmmathbf{a} \in \mathcal{A}} |N_{\mathcal{B}}^{\lambda_1}
   (\tmmathbf{a}) | |N_{\mathcal{C}}^{\lambda_2} (\tmmathbf{a}) | -
   \frac{|\mathcal{B}|}{q} \sum_{\tmmathbf{a} \in \mathcal{A}}
   |N_{\mathcal{C}}^{\lambda_2} (\tmmathbf{a}) | - \frac{|\mathcal{C}|}{q}
   \sum_{\tmmathbf{a} \in \mathcal{A}} |N_{\mathcal{B}}^{\lambda_1}
   (\tmmathbf{a}) | + \frac{|\mathcal{A}| |\mathcal{B}| |\mathcal{C}|}{q^2}
   \right| < q^{d - 1} \sqrt{|\mathcal{B}| |\mathcal{C}|} . \end{equation}
From Theorem \ref{12}, we have
\begin{eqnarray}
  \left| \sum_{\tmmathbf{a} \in \mathcal{A}} |N_{\mathcal{C}}^{\lambda_2}
  (\tmmathbf{a}) | - \frac{|\mathcal{A}| |\mathcal{C}|}{q} \right| & \leq
  & \sqrt{q^{d - 1} |\mathcal{A}| |\mathcal{C}|} \label{v2}\\
  \left| \sum_{\tmmathbf{a} \in \mathcal{A}} |N_{\mathcal{B}}^{\lambda_1}
  (\tmmathbf{a}) | - \frac{|\mathcal{A}| |\mathcal{B}|}{q} \right| & \leq
  & \sqrt{q^{d - 1} |\mathcal{A}| |\mathcal{B}|}. \label{v3}
\end{eqnarray}
Putting (\ref{v1}), (\ref{v2}) and (\ref{v3}) together, it follows that
\[ \left| \sum_{\tmmathbf{a} \in \mathcal{A}} |N_{\mathcal{B}}^{\lambda_1}
   (\tmmathbf{a}) | |N_{\mathcal{C}}^{\lambda_2} (\tmmathbf{a}) | -
   \frac{|\mathcal{A}| |\mathcal{B}| |\mathcal{C}|}{q^2} \right| \leq
   \frac{|\mathcal{B}|}{q} \sqrt{q^{d - 1} |\mathcal{A}| |\mathcal{C}|} +
   \frac{|\mathcal{C}|}{q} \sqrt{q^{d - 1} |\mathcal{A}| |\mathcal{B}|} + q^{d
   - 1} \sqrt{|\mathcal{B}| |\mathcal{C}|}, \]
completing the proof of the theorem.

\section{The system of three equations and three variables}

\subsection{The case $d = 2$ (Proof of Theorem \ref{332})}

Let $\mathcal{A}^{\ast} =\mathcal{A} \cap \mathbb{F}_q^{\ast} \times
\mathbb{F}_q^{\ast}$, $\mathcal{B}^{\ast} =\mathcal{B} \cap
\mathbb{F}_q^{\ast} \times \mathbb{F}_q^{\ast}$ and $\mathcal{C}^{\ast}
=\mathcal{C} \cap \mathbb{F}_q^{\ast} \times \mathbb{F}_q^{\ast}$ then
\[ |\mathcal{A}^{\ast} |, |\mathcal{B}^{\ast} |, |\mathcal{C}^{\ast} | \gg
   q^{3 / 2} . \]
For any $\lambda_1, \lambda_2 \in \mathbb{F}_q^{\ast}$, it follows from Theorem \ref{23}
that
\[ |\{(\tmmathbf{a}, \tmmathbf{b}, \tmmathbf{c}) \in \mathcal{A}^{\ast} \times
   \mathcal{B}^{\ast} \times \mathcal{C}^{\ast} : B (\tmmathbf{a}, \tmmathbf{b})
   = \lambda_1, B (\tmmathbf{a}, \tmmathbf{c}) = \lambda_2 \}| = (1 + o (1))
   \frac{|\mathcal{A}^{\ast} | |\mathcal{B}^{\ast} | |\mathcal{C}^{\ast}
   |}{q^2} . \]
By the pigeon-hole principle, there exists $\tmmathbf{a}_0 \in
\mathcal{A}^{\ast}$ such that
\[ |\{(\tmmathbf{b}, \tmmathbf{c}) \in \mathcal{B}^{\ast} \times
   \mathcal{C}^{\ast} : B (\tmmathbf{a}_0, \tmmathbf{b}) = \lambda_1, B
   (\tmmathbf{a}_0, \tmmathbf{c}) = \lambda_2 \}| = (1 + o (1))
   \frac{|\mathcal{B}^{\ast} | |\mathcal{C}^{\ast} |}{q^2} \gg q. \]
Let $\delta = \sqrt{|\mathcal{B}^{\ast} | |\mathcal{C}^{\ast} |} / q^2 \gg
q^{- 1 / 2}$. Let 
\[ \mathcal{B}_1 =\{\tmmathbf{b} \in \mathcal{B}^{\ast} : B (\tmmathbf{a}_0,
\tmmathbf{b}) = \lambda_1 \}, \,\,\, \mathcal{C}_1 =\{\tmmathbf{c} \in
\mathcal{C}^{\ast} : B (\tmmathbf{a}_0, \tmmathbf{c}) = \lambda_2 \}, \] then
$|\mathcal{B}_1 | |\mathcal{C}_1 | \gg \delta^2 q^2$. We assume that
$|\mathcal{C}_1 | \geq |\mathcal{B}_1 |$, then $|\mathcal{C}_1 |
\geq \delta q$. It suffices to show that there are at least $\delta q$
values of $\lambda$ such that
\[ |\{(\tmmathbf{b}, \tmmathbf{c}) \in \mathcal{B}_1 \times \mathcal{C}_1 : B
   (\tmmathbf{b}, \tmmathbf{c}) = \lambda\}| > 0. \]
   For a fix $\tmmathbf{b}_0 \in \mathcal{B}_1$, we want to
solve the following system
\begin{equation}\label{t1} B (\tmmathbf{a}_0, \tmmathbf{c}) = \lambda_2, B (\tmmathbf{b}_0, \tmmathbf{c})
   = \lambda, \tmmathbf{c} \in \mathcal{C}_1, \end{equation}
under the constraint $B (\tmmathbf{a}_0, \tmmathbf{b}_0) = \lambda_1$. Suppose
that $B (\tmmathbf{x}, \tmmathbf{y}) = x_1 y_1 + \kappa x_2 y_2$ for some
$\kappa \in \mathbb{F}_q^{\ast}$. Let $\tmmathbf{a}_0 = (a_1, a_2)$,
$\tmmathbf{b}_0 = (b_1, b_2)$ and $\tmmathbf{c}= (c_1, c_2)$, the system (\ref{t1})
becomes
\begin{eqnarray*}
  a_1 b_1 + \kappa a_2 b_2 & = & \lambda_1\\
  a_1 c_1 + \kappa a_2 c_2 & = & \lambda_2\\
  b_1 c_1 + \kappa b_2 c_2 & = & \lambda .
\end{eqnarray*}
This implies that $\kappa (a_2 b_1 - a_1 b_2) c_2 = \lambda_2 b_1 - \lambda
a_1$. Thus the system has at most one solution if $a_2 b_1 - a_1 b_2 \neq 0$.
If $a_2 b_1 = a_1 b_2$, the system is solvable only if $\lambda =
\lambda_2 b_1 / a_1$. Besides, from the first equation, $(\kappa a_2^2 +
a_1^2) b_2 = \lambda_1 a_2$ and $(\kappa a_2^2 + a_1^2) b_1 = \lambda_1 a_1$.
It follows that the system (\ref{t1}) has at most one solution if
\begin{equation}\label{t2} \tmmathbf{b}_0 = (b_1, b_2) \neq \left( \frac{\lambda_1 a_1}{\kappa a_2^2 +
   a_1^2}, \frac{\lambda_1 a_2}{\kappa a_2^2 + a_1^2} \right) . \end{equation}
Since $|\mathcal{C}^{*}| \leq q$, $|\mathcal{B}_1| \geq \delta^2 q \gg 1$.      
Thus, we can choose $\tmmathbf{b}_0 \in \mathcal{B}_1$
satisfying (\ref{t2}). The system (\ref{t1}) has at most one solution for each $\lambda$. So there exists at least $|\mathcal{C}_1 | \geq \delta q$ values of $\lambda$ such
that the system (\ref{t1}) is solvable. This completes the proof of the theorem.

\subsection{The case $d \geq 3$ (Proof of Theorem \ref{333})}

Let $\mathcal{A}^{\ast} =\mathcal{A} \backslash (0, \ldots, 0)$,
$\mathcal{B}^{\ast} =\mathcal{B} \backslash (0, \ldots, 0)$ and
$\mathcal{C}^{\ast} =\mathcal{C} \backslash (0, \ldots, 0)$. For any
$\lambda_1, \lambda_2 \in \mathbb{F}_q^{\ast}$, it follows from Theorem \ref{23} that
\[ |\{(\tmmathbf{a}, \tmmathbf{b}, \tmmathbf{c}) \in \mathcal{A}^{\ast} \times
   \mathcal{B}^{\ast} \times \mathcal{C}^{\ast} : B (\tmmathbf{a}, \tmmathbf{b})
   = \lambda_1, B (\tmmathbf{a}, \tmmathbf{c}) = \lambda_2 \}| = (1 + o (1))
   \frac{|\mathcal{A}^{\ast} | |\mathcal{B}^{\ast} | |\mathcal{C}^{\ast}
   |}{q^2} . \]
By the pigeon-hole principle, there exists $\tmmathbf{a}_0 \in
\mathcal{A}^{\ast}$ such that
\[ |\{(\tmmathbf{b}, \tmmathbf{c}) \in \mathcal{B}^{\ast} \times
   \mathcal{C}^{\ast} : B (\tmmathbf{a}_0, \tmmathbf{b}) = \lambda_1, B
   (\tmmathbf{a}_0, \tmmathbf{c}) = \lambda_2 \}| = (1 + o (1))
   \frac{|\mathcal{B}^{\ast} | |\mathcal{C}^{\ast} |}{q^2} \gg q^d . \]
For any $\tmmathbf{a} \in \mathbb{F}_q^d \backslash (0, \ldots, 0)$, set
$\Pi_{\lambda} (\tmmathbf{a}) =\{\tmmathbf{v} \in \mathbb{F}_q^d : B
(\tmmathbf{a}, \tmmathbf{v}) = \lambda\}$. Let 
\[ \mathcal{B}_1 = \Pi_{\lambda_1}
(\tmmathbf{a}_0) \cap \mathcal{B}^{\ast}, \,\,\, \mathcal{C}_1 = \Pi_{\lambda_2}
(\tmmathbf{a}_0) \cap \mathcal{C}^{\ast},\] then $|\mathcal{B}_1 | |\mathcal{C}_1
| \gg q^d$. Theorem \ref{333} follows immediately from the following lemma.

\begin{lemma}
  For any $\tmmathbf{a} \in \mathbb{F}_q^d \backslash (0, \ldots, 0)$ and
  $\lambda_1, \lambda_2 \in \mathbb{F}_q^{\ast}$, suppose that $\mathcal{E}
  \subseteq \Pi_{\lambda_1} (\tmmathbf{a})$, $\mathcal{F} \subseteq
  \Pi_{\lambda_2} (\tmmathbf{a})$. If $d \geq 3$ and $|\mathcal{E}|
  |\mathcal{F}| \gg q^d$, then
  \[ | \Pi (\mathcal{E}, \mathcal{F}) : =\{B (\tmmathbf{e}, \tmmathbf{f}) :
     \tmmathbf{e} \in \mathcal{E}, \tmmathbf{f} \in \mathcal{F}\}| \geq (1
     - o (1)) q. \]
\end{lemma}

\begin{proof} The proof is similar to that of \cite[Theorem 2.8]{hart}. Define the incidence function
  \[ v_{\lambda} (\mathcal{E}, \mathcal{F}) =\{(\tmmathbf{e}, \tmmathbf{f}) \in
     \mathcal{E} \times \mathcal{F}: B (\tmmathbf{e}, \tmmathbf{f}) = \lambda\}.
  \]
  The Fourier transform of a complex-valued function $f$ on $\mathbb{F}_q^d$
with respect to a non-trivial additive character $\chi$ on
$\mathbb{F}_q$ is given by
\[ \hat{f} (k) = q^{- d} \sum_{x \in \mathbb{F}_q^d} \chi (- x \cdot k) f
   (x), \]
and the Fourier inversion formula takes the form
\[ f (x) = \sum_{k \in \mathbb{F}_q^d} \chi (x \cdot k) \hat{f} (k) . \]
  Using exponential sums, Hart, Iosevich, Koh, and Rudnev (\cite[Theorem 2.1]{hart}) showed that
  \[ \sum_{\lambda \in \mathbb{F}_q} v_{\lambda}^2 (\mathcal{E}, \mathcal{F})
     \leq |\mathcal{E}|^2 |\mathcal{F}|^2 q^{- 1} + |\mathcal{E}|q^{2 d -
     1} \sum_{\tmmathbf{f} \in \mathbb{F}_q^d \backslash (0, \ldots 0)}
     |\mathcal{F} \cap l_{\tmmathbf{f}} | | \hat{\mathcal{F}} (\tmmathbf{f}) |^2
     + (q - 1) q^{- 1} |\mathcal{E}| |\mathcal{F}|\mathcal{F}((0, \ldots, 0)),
  \]
  where
  \[ l_{\tmmathbf{f}} =\{t\tmmathbf{f}: t \in \mathbb{F}_q^{\ast} \}. \]
  Since $\mathcal{F} \subseteq \Pi_{\lambda_2} (\tmmathbf{a})$, $(0, \ldots, 0)
  \notin \mathcal{F}$ and $|\mathcal{F} \cap l_{\tmmathbf{f}} | \leq 1$ for
  any $\tmmathbf{f} \in \mathbb{F}_q^d \backslash (0, \ldots 0)$. Therefore,
  \begin{eqnarray*}
    \sum_{\lambda \in \mathbb{F}_q} v_{\lambda} (\mathcal{E}, \mathcal{F})^2
    & \leq & |\mathcal{E}|^2 |\mathcal{F}|^2 q^{- 1} + |\mathcal{E}|q^{2
    d - 1} \sum_{\tmmathbf{f} \in \mathbb{F}_q^d \backslash (0, \ldots 0)} |
    \hat{\mathcal{F}} (\tmmathbf{f}) |^2\\
    & \leq & |\mathcal{E}|^2 |\mathcal{F}|^2 q^{- 1} + |\mathcal{E}|q^{2
    d - 1} q^{- d} \sum_{\tmmathbf{f}' \in \mathbb{F}_q^d}
    \mathcal{F}(\tmmathbf{f}')^2\\
    & = & |\mathcal{E}|^2 |\mathcal{F}|^2 q^{- 1} + |\mathcal{E}|
    |\mathcal{F}|q^{d - 1} .
  \end{eqnarray*}
  By the Cauchy-Schwartz inequality, we have
  \[ |\mathcal{E}|^2 |\mathcal{F}|^2 = \left( \sum_{\lambda} v_{\lambda}
     (\mathcal{E}, \mathcal{F}) \right)^2 \leq | \Pi (\mathcal{E},
     \mathcal{F}) | \sum_{\lambda} v_{\lambda} (\mathcal{E}, \mathcal{F})^2 .
  \]
  This implies that
  \[ | \Pi (\mathcal{E}, \mathcal{F}) | \geq \frac{q}{1 +
     \frac{q^d}{|\mathcal{E}| |\mathcal{F}|}} . \]
  This follows that if $|\mathcal{E}| |\mathcal{F}| \gg q^d$ then $| \Pi
  (\mathcal{E}, \mathcal{F}) | = q (1 - o (1))$, completing the proof of the
  lemma.
\end{proof}

\bibliographystyle{amsplain}

\begin{thebibliography}{10}

\bibitem{cipra} J. A. Cipra, T. Cochrane and C. Piner, \textit{Heilbronn's conjecture on Waring's number (mod $p$)}, J. Number Theory \textbf{125}(2) (2007), 289--297.

\bibitem{covert} D. Covert, D. Hart, A. Iosevich, and I. Uriarte-Tuero, \textit{An analog of the Furstenberg-Katznelson-Weiss theorem on traingles in sets of positive density in finite field geometries}, preprint 2008, arXiv:0804.4894.

\bibitem{s2} K. Gyarmati and A. S\'ark\"ozy, \textit{Equations in finite fields with restricted solution sets, II (algebraic equations)}, Acta Math. Hungar. \textbf{119} (2008), 259–-280.

\bibitem{hart-thesis} D. Hart, \textit{Explorations of Geometric Combinatorics in Vector Spaces over Finite Fields}, PhD Thesis, Missouri University.

\bibitem{hart-iosevich} D. Hart and A. Iosevich, \textit{Sums and products in finite fields: an integral geometric viewpoint}, Contemp. Math. \textbf{464} (2008).

\bibitem{hart2} D. Hart and A. Iosevich, \textit{Ubiquity of simplices in vector spaces over finite fields}, Anal. Math. \textbf{34}(1) (2008).

\bibitem{hart} D. Hart, A. Iosevich, D. Koh and M. Rudnev, \textit{Averages over hyperplanes, sum-product theory in finite fields, and the Erd\H{o}s-Falconer distance conjecture}, to appear in Trans. Amer. Math. Soc., arXiv:0707.3473.

\bibitem{hart3} D. Hart, A. Iosevich, D. Koh, S. Senger, and I. Uriarte-Tuero, \textit{Distance graphs in vector spaces over finite fields, coloring, pseudo-randomness and arithmetic progressions}, preprint 2008, arXiv:0804.3036.

\bibitem{pseudo} M. Krivelevich and B. Sudakov, \textit{Pseudo-random graphs}, in More Sets, Graphs and Numbers, Bolyai Soc. Math. Studies 15, Springer, 2006, 199-262. 

\bibitem{igor} I. E. Shparlinski, \textit{On the solvability of bilinear equations in finite fields}, Glasg. Math. J. \textbf{50} (2008), 523--529.

\bibitem{s1} A. S\'ark\"ozy, \textit{On products and shifted products of residues modulo $p$,} INTEGERS, \textbf{8}(2) (2008), A9.

\bibitem{vinh1} L. A. Vinh, \textit{On a Furstenberg-Katznelson-Weiss type theorem over finite fields}, to appear in Ann. Comb., arXiv:0807.2849

\bibitem{vinh2} L. A. Vinh, \textsl{On kaleidoscopic pseudo-randomness of finite Euclidean graphs}, preprint 2008, arXiv:0807.2689. 

\bibitem{vinh3} L. A. Vinh, \textit{On $k$-simplexes in $(2k-1)$-dimensional vector spaces over finite fields}, to appear in Proc. 21st FPSAC 2009.

\bibitem{vinh4} L. A. Vinh, \textit{Triangles in vector spaces over finite fields}, to appear in Online J. Anal. Comb. (2009).

\bibitem{weil} A. Weil, \textit{Number of solutions of equations in finite fields}, Bull. AMS \textbf{55} (1949), 497--508.

\end{thebibliography}

\end{document}